\documentclass[12pt,titlepage]{utarticle}
\usepackage{amsmath,amssymb,amsthm,pictexwd,graphicx,color}
\usepackage{setspace}
\usepackage[mathscr]{eucal}
\usepackage[colorlinks=true, pdfstartview=FitV, linkcolor=blue, citecolor=blue, urlcolor=blue]{hyperref}

\numberwithin{equation}{section}

\def\BZ{\mathbb{Z}}

\newcommand{\MR}[1]{{\mathbb{R}^{#1}}}            % Real numbers
\newtheorem{defn}{Definition}

\newtheorem{lma}{Lemma}
\newtheorem{ilma}{Lemma}
\newtheorem{prop}{Proposition}
\newtheorem{iprop}{Proposition}
\newtheorem{corr}{Corollary}

\def\Spec{\mathrm{Spec}}

\def\Hom{\mathrm{Hom}}
\def\End{\mathrm{End}}

\def\Ext{\mathrm{Ext}}
\def\Coh{\mathrm{Coh}}
\def\QCoh{\mathrm{QCoh}}

\def\dim{\mathrm{dim}}
\def\mcO{\mathcal{O}}
\def\mcA{\mathcal{A}}
\def\mcB{\mathcal{B}}
\def\mcC{\mathcal{C}}
\def\mcD{\mathcal{D}}
\def\mcE{\mathcal{E}}
\def\mcF{\mathcal{F}}

\def\AMod{{A-\mathrm{Mod}}}
\def\Amod{{A-\mathrm{mod}}}

\def\ie{\textit{i.e.,\ }}

\begin{document}
\preprint{MIFP--08--05\\}

\title{A Note on Support in Triangulated Categories}

\author{Aaron Bergman\address{George P. \& Cynthia W. Mitchell Institute for Fundamental Physics\\
        Texas A\&M University\\
        College Station, TX 77843-4242\\ {~}\\
        \email{abergman@physics.tamu.edu}}}

\Abstract{In this note, I define a notion of a compactly supported object in a triangulated category. I prove a number of propositions relating this to traditional notions of support and give an application to the theory of derived Morita equivalence. I also discuss a connection to supersymmetric gauge theories arising from D-branes at a singularity.}

\maketitle
\newpage

\section{Introduction}\label{sec:intro}
\subsection{Overview}

Let $\mcC$ be a triangulated category with small coproducts. We will always work over a field $k$, and all dimensions will be as vector spaces over that field. Recall the following definition:

\begin{defn}
An object $\mcE \in \mcC$ is called \textbf{compact} if $\Hom(\mcE,-)$ commutes with coproducts.
\end{defn}

The full subcategory of compact objects is triangulated, and we will denote it as $\mcC^c$. If $\mcC \cong \mcD(\QCoh(X))$ for some variety $X$, then it is well-known (see, for example, \cite{Neeman:1996gd}) that the compact objects are those complexes of sheaves locally quasi-isomorphic to bounded complexes of locally free coherent sheaves\footnote{If X is noetherian, $\Coh(X)$ is an abelian category, and $\mcD(\Coh(X))$ is equivalent to the full subcategory of $\mcD(\mathrm{QCoh}(X))$ whose cohomology sheaves are coherent. It is straightforward to see that the cohomology sheaves of any compact complex are coherent, so the compact objects can be considered as a subcategory of $\mcD(\Coh(X))$.}. If, on the other hand, $\mcC \cong \mcD(\AMod)$ for an algebra $A$, then the compact objects are those complexes of modules quasi-isomorphic to bounded complexes of finitely generated projective modules\footnote{Here $\AMod$ is the abelian category of possibly infinitely generated modules. If $A$ is noetherian, finitely generated modules form an abelian category which we will denote $\Amod$, and a remark similar to the previous footnote applies.}   \cite{Rickard:1988mt}.  In both cases, these are often termed \textit{perfect} objects. We will denote the triangulated subcategory of compact objects as $\mcD^c(\AMod)$ or $\mcD^c(\QCoh(X))$. We will also sometimes refer to an object in an abelian category as compact or perfect if its image in the derived category is compact.

In this note, we will consider a notion of support for an object in a triangulated category. In particular, we define

\begin{defn}
An object $\mcE \in \mcC$ is said to be \textbf{compactly supported} if, for all compact $\mcF \in \mcC$,
$$
\sum_{i=-\infty}^\infty \mathrm{dim\ }\Hom(\mcF,\mcE[i]) < \infty \ .
$$
\end{defn}

This is equivalent to saying that the homological functor $\mcC(-,\mcE)$ is of finite type when restricted to the compact subcategory. Given a triangulated category, we can consider the full subcategory of compactly supported objects which we will denote by $\mcC_{cs}$. For our examples, we will denote these subcategories as $\mcD_{cs}(\QCoh(X))$ and $\mcD_{cs}(\AMod)$.
\begin{prop}
\label{fullsub}
The full subcategory of compactly supported objects is a triangulated subcategory.
\end{prop}
\begin{proof}
It is obviously closed under shifts and finite sums. That it is closed under extensions can be seen from the long exact sequence in Homs arising from a distinguished triangle.
\end{proof}

From the point of view of noncommutative geometry, the more relevant object is often the compact subcategory of a given triangulated category. Thus, it is interesting to consider the resulting subcategories of compactly supported compact objects. Then, the above definition reduces to the statement that $\mcC^c(-,\mcE)$ is of finite type. If $\mcC^c_{cs}$ has a generator, then it is a compact algebraic $k$-linear space in the sense of Kontsevich\footnote{More properly, in order to be a linear space, $\mcC$ should be enhanced over dg-Vect, and we can apply our criterion to the dimensions of the cohomologies of the Hom-complexes.} \cite{Kontsevich:2008sl}.

To relate this to more usual senses of support, we have:
\newcounter{supsheafc}
\setcounter{supsheafc}{\value{prop}}
\begin{prop}
\label{supsheaf}
Let $X$ be a noetherian variety and $\mcE$ a bounded complex of coherent sheaves such that, for all $i\in\BZ$, the support of $H^i(\mcE)$ is proper. Then $\mcE$ is compactly supported in $\mcD(\Coh(X))$.
\end{prop}

\newcounter{supmodc}
\setcounter{supmodc}{\value{prop}}
\begin{prop}
\label{supmod}
Let $A$ be an algebra over a field, and $\mcE$ a bounded complex of modules. If
$$
\dim\ H^i(\mcE) < \infty\quad\mathrm{for\ all\ }i\in\BZ\ ,
$$
then $\mcE$ is compactly supported in $\mcD(\AMod)$.
\end{prop}

\noindent These two propositions will be proven in section \ref{sec:supp}, along with a characterization of compact support when $\mcC$ has a compact generator. To avoid confusion (or perhaps to foster it), we will always use `proper support' to refer to the support of a sheaf being a proper subscheme, and `compact support' to refer to the homological notion defined above.

With more restrictive assumptions, we can prove converses of the above propositions.

\newcounter{modulec}
\setcounter{modulec}{\value{prop}}
\begin{prop}
\label{module}
Let $A$ be an algebra over a field. Then $\mcD^c_{cs}(\AMod)$ is equivalent to the full subcategory of $\mcD^c(\AMod)$ consisting of objects whose cohomology modules have finite dimension.
\end{prop}
\newcounter{sheafc}
\setcounter{sheafc}{\value{prop}}
\begin{prop}
\label{sheaf}
Let $X$ be a smooth noetherian variety. Then $\mcD^c_{cs}(\Coh(X))$ is equivalent to the full subcategory of $\mcD^c(\Coh(X))$ consisting of objects whose cohomology sheaves have proper support.
\end{prop}

For both propositions, we can reduce to the case of specific objects in the abelian category through the use of the hypercohomology spectral sequence and repeated truncations. For the latter proposition, we also require the following lemma:

\newcounter{suppc}
\setcounter{suppc}{\value{lma}}
\begin{lma}
\label{supp}
Let $X$ be a smooth noetherian variety. Then a coherent sheaf $\mcE$ has proper support if and only if, for all coherent sheaves $\mcF$,
$$
\mathrm{dim\ }\Ext_X^i(\mcF,\mcE) < \infty \quad\mathrm{for\ all\ }i\in\BN\ .
$$
\end{lma}
\noindent This will follow from applying the hypothesis to the structure sheaves of reduced curves in the support of $\mcE$ and applying the valuative criteria of properness.

\subsection{Physics motivation}

The questions addressed in this note arise in the study of D-branes at singularities in string theory. In particular, let $X$ be a del Pezzo surface and $K_X$ be the total space of the canonical bundle on $X$. $K_X$ has a trivial canonical bundle and is thus a non-compact Calabi-Yau. It serves as a local model for a singularity as one can collapse the zero section to obtain a singular cone. In string theory, this collapse can be interpreted as a deformation of the K\"ahler metric and is, as such, not visible to the topological B-model string. Thus, to make sense of the topological B-model on the singular space, it suffices to study it on this smooth resolution.

As has been conjectured in \cite{Douglas:2000gi,Kontsevich:1994ho,Sharpe:1999qz}, the topological B-model is described by the (A$_\infty$ enhancement of the) bounded derived category of coherent sheaves\footnote{When our target is smooth, this is equivalent to the perfect subcategory of the derived category.}. Using the resulting dictionary, we can turn many physical results into mathematical theorems and vice versa. For example, the independence of the topological B-model of the choice of resolution of the singularity becomes a statement about equivalences of derived categories. In particular, Bridgeland \cite{Bridgeland:2002fd} has shown that the derived categories of coherent sheaves of different crepant resolutions of a projective three-fold with terminal singularities are equivalent.

To study D-branes located at a singularity, we need to study the stability of skyscraper sheaves located on the zero section of $K_X$. This is done in \cite{Bergman:2007xb} following results of Bridgeland \cite{Bridgeland:2006sc}. In particular, we will let $T$ be a locally free generator\footnote{Such generators may be obtained from pulling back strong, simple exceptional collections on $X$. See, for example, \cite{BridgeTStruct,Bergman:2005mz}.} of $\mcD(\Coh(K_X))$ such that $\Ext^i(T,T) = 0$ for $i\neq 0$. Let $A = \End(T)^\mathrm{op}$. Then Rickard's derived Morita equivalence \cite{Rickard:1988mt} tells us that there is an equivalence of categories:
\begin{equation}
\label{equiv}
\mcD^b(\Coh(K_X)) \cong \mcD^b(\Amod)
\end{equation}
where the latter category is the bounded derived category of finitely generated $A$-modules. Here, the bounded derived categories arise as the perfect subcategories considered above as $K_X$ is smooth noetherian and $A$ is noetherian of finite global dimension. $A$ is a Calabi-Yau algebra \cite{Bocklandt:2006gc,Ginzburg:2006cy} and can be represented as the path algebra of a quiver with relations derived from a potential. This data, along with the dimension vector obtained from the module corresponding to a skyscraper sheaf, defines the holomorphic data of an $\mathcal{N} = 1$ quiver gauge theory. This gauge theory is the worldvolume theory that lives on a D3-brane located at the ``tip" of the cone, \ie the D3-brane is given by the embedding of $\MR{4}$ into to $\MR{4} \times \mcC(X)$ where $\mcC(X)$ is the cone obtained by collapsing the zero section of $K_X$, and we embed at the tip of $\mcC(X)$.

It is an interesting question to ask what is the proper ``physical'' category to describe the topological B-model on these noncompact spaces. The category $\mcD^b(\Coh(K_X))$ is unsuitable because, among other things, it is not Calabi-Yau. In particular, there is no Serre functor as $K_X$ is noncompact. In the mathematics literature, it is common to study the full subcategory of objects whose cohomology sheaves have reduced support along the zero section of $K_X$. Under the equivalence \eqref{equiv}, this maps to the full subcategory of nilpotent modules \cite{BridgeTStruct}. This is a Calabi-Yau category and has many nice mathematical properties, but its physical significance is still somewhat mysterious. Another obvious category to consider is $\mcD^b_{\mathrm{prop}}(\Coh(K_X))$, the full subcategory with proper supports for the cohomology sheaves. This is also Calabi-Yau in the sense that it has a Serre functor equivalent to the shift by three functor. It was conjectured in \cite{Bergman:2007xb} that this is equivalent to the category $\mcD^b_\mathrm{fd}(\AMod)$ consisting of objects whose cohomology modules have finite dimension vector. Using the above propositions, we can characterize these subcategories intrinsically as the full subcategories with compact support. Thus, we have

\begin{corr}
Let $X$ be a smooth, noetherian variety, and let $T$ be a compact generator of $\mcD(Coh(X))$ such that $\Hom(T,T[i])=0$ for $i\neq0$ and $End(T)$ is left-noetherian and of finite global dimension. Then
$$
\mcD^b_\mathrm{prop}(Coh(X)) \cong \mcD^b_\mathrm{fd}(\End(T)^\mathrm{op}\mathrm{-mod})\ .
$$
\end{corr}

\subsection{Acknowledgments}

I would like to thank David Ben-Zvi for his help with the manuscript and for many other suggestions during the writing of this note. I would also like to thank Zach Teitler for a suggestion that led to a proof of lemma \ref{supp} and Tom Nevins for helpful e-mails. This research was supported by the National Science Foundation under Grant No. PHY-0505757, and by Texas A\&M University.

\section{Proofs of propositions \ref{supsheaf} and \ref{supmod}}\label{sec:supp}

For convenience, we restate proposition \ref{supsheaf}:

\setcounter{iprop}{\value{supsheafc}}
\begin{iprop}
Let $X$ be a noetherian variety and $\mcE$ a bounded complex of coherent sheaves such that, for all $i\in\BZ$, the support of $H^i(\mcE)$ is proper. Then $\mcE$ is compactly supported in $\mcD(\Coh(X))$.
\end{iprop}
\begin{proof}
We begin with the case where $\mcE$ is a sheaf. We can consider $\mcE$ as $p_*\mcE$ where $p : \mathrm{Supp\ }\mcE \to X$ is the embedding of the support with the natural subscheme structure. Then, for a perfect $\mcF$ we can compute
$$
\Hom_{\mcD^b(X)}(\mcF,\mcE[i]) \cong \Hom_{\mcD^b(\Coh(\mathrm{Supp\ }\mcE))}(\mathbf{L}p^*\mcF,\mcE[i])
\cong H^i_{\mathrm{Supp\ }\mcE}(\mathbf{R}\mathcal{H}om(\mathbf{L}p^*\mcF,\mcE))\ .
$$
Derived pullbacks preserve perfectness (cf.\ \cite{TT} 2.5.1). As $X$ is finite type, any perfect complex is bounded, and local Homs from perfect complexes preserve boundedness, so $\mathbf{R}\mathcal{H}om(\mathbf{L}p^*\mcF,\mcE)$ is a bounded complex of sheaves. Since $\mathrm{Supp\ }\mcE$ is finite dimensional and noetherian, there are only a finite number of non-zero terms in the hypercohomology spectral sequence, and they are all finite dimensional as $\mathrm{Supp\ }\mcE$ is proper. Thus, $\mcE$ is compactly supported. For $\mcE$ a bounded complex, the result follows from examining the analog of the hypercohomology spectral sequence for the hyperext functor $\Ext^i(\mcF_\bullet,-)$ and applying it to $\mcE$, reducing to the previous case applied to the finite number of non-zero cohomology sheaves.
\end{proof}

The proof of proposition \ref{supmod} is similar:
\setcounter{iprop}{\value{supmodc}}
\begin{iprop}
Let $A$ be an algebra over a field, and $\mcE$ a bounded complex of modules. If
$$
\dim\ H^i(\mcE) < \infty\quad\mathrm{for\ all\ }i\in\BZ\ ,
$$
then $\mcE$ is compactly supported in $\mcD(\AMod)$.
\end{iprop}
\begin{proof}
We again start with the case of a finite dimensional module, $M$. We can write its dimension as $\dim\ \Hom_{\AMod}(A,M)$. Since any finitely generated projective module is a direct summand in a free module on a finite set, we immediately have that
$$
\dim\ \Hom_\AMod(P,M) < \infty
$$
for all finitely generated projective modules, $P$. Since a perfect complex of modules is quasiisomorphic to a bounded complex of finitely generated projective modules, we see that any finite dimensional module is compactly supported. We can again extend this to bounded complexes by usual spectral sequence arguments.
\end{proof}

In the case where a triangulated category has a compact generator there is useful characterization of compact support. As above, we let $\mcC$ denote our triangulated category. An object $\mcA$ is called a compact generator if it is compact and if $\Hom_\mcC(\mcA,\mcE[i]) = 0$ for all $i$ implies that $\mcE\cong0$. A theorem of Ravenel and Neeman \cite{Neeman:1992kl} implies that the smallest triangulated subcategory of $\mcC$ containing $\mcA$ and which is closed under direct summands is $\mcC^c$. In this situation, there is the following criterion for compact support:
\begin{lma}
\label{genlma}
An object, $\mcE$, in a triangulated category, $\mcC$, with compact generator, $\mcA$, is compactly supported if and only if
$$
\sum_{i=-\infty}^{\infty}\dim\ \Hom_\mcC(\mcA,\mcE[i]) < \infty\ .
$$
\end{lma}
\begin{proof}
One direction is obvious. To see the other direction, let $\mcB$ be the full subcategory of objects in $\mcF\in\mcC$ such that
$$
\sum_{i=-\infty}^{\infty}\dim\ \Hom_\mcC(\mcF,\mcE[i]) < \infty\ .
$$
By assumption, $\mcA$ is an object in $\mcB$. Furthermore, as in proposition \ref{fullsub}, it follows from the long exact sequence associated to $\Hom_\mcC(-,\mcE[i])$ that $\mcB$ is a triangulated subcategory. Finally, it is obviously closed under direct summands. Thus, by Ravenel and Neeman's theorem, $\mcB$ must contain $\mcC^c$, and $\mcE$ is compactly supported.
\end{proof}

We can now give an essentially equivalent proof of proposition \ref{supmod} by noting that $\mcD(\AMod)$ is generated by the free module $A$ (for a longer discussion, see, for example, \cite{Keller:1998te}) and computing
\begin{equation}
\label{genrel}
\sum_{i=-\infty}^{\infty}\dim\ \Hom(A,\mcE[i]) = \sum_{i=-\infty}^{\infty}\dim\ H^i(\mcE) < \infty\ .
\end{equation}
In fact, this suffices to also prove proposition \ref{module}, but we prove it by another method in the following section that generalizes to a proof of proposition \ref{sheaf}.

In addition, a theorem of Bondal and van den Bergh \cite{Bondal:2003gr} states that, for any quasi-compact separated scheme, X, $\mcD(\QCoh(X))$ is generated by a single perfect complex, giving an often simple criterion for compact support.

\section{Proof of proposition \ref{module}}

%We first note the following obvious consequence of the hypercohomology spectral sequence:

%\begin{lma}
%\label{hyperco}
%Let $\mcA$ be an abelian category and $\mcD(\mcA)$ be its derived category. Then, for $\mcE \in \mcD^b(\mcA)$ and $\mcF \in \mcA$, if
%$$
%\sum_{j=-\infty}^\infty \sum_{i=0}^\infty \mathrm{dim\ }\Ext_\mcA^i(\mcF,H^j(\mcE)) < \infty\ .
%$$
%then
%$$
%\sum_{i=-\infty}^\infty \mathrm{dim\ }\Hom_{\mcD(\mcA)}(\mcF,\mcE[i]) < \infty
%$$
%\end{lma}

%\begin{proof} There is a spectral sequence converging to $\Hom_{\mcD(\mcA)}(\mcF,\mcE[i])$ with
%$$
%E_2^{pq} = \Ext_\mcA^p(\mcF,H^{q}(\mcE))\ .
%$$
%Then, by assumption, these groups are all finite, and they are nonzero for only a finite number of $p,q$. Thus, we have a finite sum of finite terms, and we are done.
%\end{proof}

We begin with the following lemma:

\begin{lma}
\label{conproj}
Let $A$ be an algebra over a field, $P$ a projective module and $\mcE$ a bounded complex of modules. Then
$$
\sum_{i=-\infty}^\infty \mathrm{dim\ }\Hom_{\mcD^b(\AMod)}(P,\mcE[i]) < \infty
$$
implies that
$$
\sum_{i=-\infty}^\infty \mathrm{dim\ }\Hom_\AMod(P,H^i(\mcE)) < \infty\ .
$$
\end{lma}
\begin{proof}
Because the assumptions of the lemma are invariant under shifts and $\mcE$ is bounded, we can assume that $\tau_{\ge0}\mcE \cong \mcE$ where $\tau$ is the usual truncation functor. There exists a distinguished triangle
$$
\tau_{\le 0}\mcE \longrightarrow \mcE  \longrightarrow \tau_{\ge 1} \mcE \ .
$$
By assumption $\tau_{\le 0}\mcE  \cong \tau_{\le0}\tau_{\ge0}\mcE  \cong H^0(\mcE )$. We can apply the cohomology functor $\Hom(P,-)$ to get a long exact sequence. Since $P$ is in $\mcD^{\le0}$ and all three elements in the triangle are in $\mcD^{\ge 0}$, all the negative Homs vanish. Thus, we can write
\begin{equation*}
\begin{split}
0 \longrightarrow \Hom(P,H^0(\mcE))& \longrightarrow \Hom(P,\mcE) \longrightarrow \Hom(P,\tau_{\ge 1}\mcE) \longrightarrow  \\
& \Hom(P,H^0(\mcE)[1]) \longrightarrow\Hom(P,\mcE[1]) \longrightarrow \Hom(P, (\tau_{\ge 1}\mcE)[1]) \longrightarrow \dots \ .
\end{split}
\end{equation*}

As $\Hom(P,\mcE)$ is finite dimensional, exactness implies that $\Hom(P,H^0(\mcE))$ is finite dimensional. Furthermore, $\Hom(P,H^0(\mcE)[i]) = 0$ for $i>0$ as $P$ is projective. Thus, $\Hom(P,(\tau_{\ge1}\mcE)[i])$ is also finite dimensional for all $i$. We can now repeat this argument for $\tau_{\ge1}\mcE$ and proceed inductively to see that $\Hom(P,H^i(\mcE))$ is finite dimensional for all $i$. \end{proof}

Note that the long exact sequence in fact implies that $\Hom(P,\mcE[i]) \cong \Hom(P,\tau_{\ge i}\mcE[i]) \cong \Hom(P,H^i(\mcE))$ which can also be computed directly. Thus, the sums are equal, and we reproduce the fact  \eqref{genrel} noted in the previous section.

Combining these results, we have:

\setcounter{iprop}{\value{modulec}}
\begin{iprop}
Let $A$ be an algebra over a field. Then $\mcD^c_{cs}(\AMod)$ is equivalent to the full subcategory of $\mcD^c(\AMod)$ consisting of objects whose cohomology modules have finite dimension.
\end{iprop}
\begin{proof} Given an object in $\mcD^c(\AMod)$, it is quasiisomorphic to a bounded complex of modules. If, furthermore, it has finite dimensional cohomology modules, proposition \ref{supmod} tells us that it has compact support.

Now assume that $\mcE$ is compactly supported. Then
$$
\sum_{i=-\infty}^\infty \mathrm{dim\ }\Hom(A,\mcE[i]) < \infty\ .
$$
By lemma \ref{conproj}, the sum of the dimensions of the cohomology modules is finite, and we are done.
\end{proof}

\section{Proof of proposition \ref{sheaf}}

In order to prove proposition \ref{sheaf}, we first prove lemma \ref{supp} which we restate here.

\setcounter{ilma}{\value{suppc}}
\begin{ilma}
Let $X$ be a noetherian variety. Then a coherent sheaf $\mcE$ has proper support if and only if, for all coherent sheaves $\mcF$ and $i \in\BN$,
$$
\mathrm{dim\ }\Ext_X^i(\mcF,\mcE) < \infty \ .
$$
\end{ilma}

\begin{proof}
%To prove the forward implication, first note that, as $X$ is smooth, only a finite number of terms in the sum are non-zero. Next, we can rewrite the above as
%$$
%\Hom_{\mcD(\Coh(X))}(\mcF,p_*\mcE[i])
%$$
%where we consider $\mcE$ as pushed-forward from its support which we give the  natural subscheme structure. We can then apply the adjunction to get
%$$
%\mathrm{dim\ }\Ext_X^i(\mcF,\mcE) = \dim\ \Hom_{\mcD(\Coh(\mathrm{Supp\ }\mcE))}(\mathbf{L}p^*\mcF,\mcE[i])\ .
%$$
%As $X$ is smooth, the left-derived pullback is bounded. This Hom can be computed by a spectral sequence (\cite{GelMan} Ex IV.2.2b) involving Exts between the cohomology sheaves of $\mathbf{L}p^*\mcF$ and $\mcE$, and, since $\mathrm{Supp\ }\mcE$ is proper, these are finite dimensional.

We can compute the Ext group by the local to global spectral sequence with
$$
E_2^{pq} = H^p(\mathcal{E}xt^q(\mcF,\mcE))\ .
$$
Since $\mcE$ has proper support, the local Ext also does, and thus the dimensions of the cohomology groups are finite dimensional. This is an upper-right quadrant spectral sequence, so only a finite number of terms contribute to any given $\Ext_X^i(\mcF,\mcE)$, and we are done.

For the other direction, we apply the valuative criterion for properness. Since $X$ is noetherian, $E = \mathrm{Supp\ }\mcE$ is a closed subscheme. Let $R$ be a discrete valuation ring and $K$ its field of fractions. We are given a map from $\Spec\ K$ into $E$. The closure of the image of this in $X$ map is a closed subvariety, $f :C \to E$, with the reduced induced structure. $C$ is a curve and is thus (cf. \cite{Hartshorne} Ex. IV.1.4) either proper or affine. We can restrict $\mcE$ to $C$, and we have by assumption that
$$
\dim\ H^0_C(f^!\mcE) = \dim\ \Hom_X(\mcO_C,\mcE) <\infty  \ .
$$
Since $E$ is a closed subscheme of $X$, $C$ lies in the support of $\mcE$ and thus $f^!\mcE = \mathcal{H}om(\mcO_C,\mcE)$ is globally supported on $C$. However, if $C$ were affine, this would mean that the cohomology would be infinite dimensional, so $C$ must be proper. Therefore, there exists a map from $\Spec\ R$ to $E$, and thus $E$ is proper.
%
%For the other direction, as $X$ is smooth, every sheaf is perfect, and we can let $\mcF$ be the structure sheaf, $\mcO_C$, of a reduced, rational curve contained in the support of $\mcE$. Let $p$ be a map embedding $C$ into $X$. Then $\Hom_X(p_*\mcO_C,\mcE) \cong H_C^0(p^*\mcE)$. If $C$ is non-proper, it is an affine scheme, and since it is embedded in the support of $\mcE$, $p^*\mcE$ is supported everywhere on $C$. In terms of an algebra $R$ such that $C = \Spec R$, $p^*\mcE$ corresponds to a module with trivial annihilator. However, any such module is infinite dimensional as a $\BC$-module contradicting the statement that $\dim\ \Hom_X(p_*\mcO_C,\mcE)<\infty$. Thus $C$ must be proper, meaning that the support of $\mcE$ is complete and thus is proper.
\end{proof}
\noindent Note that to prove a sheaf has proper support, we only need to test the Homs from the structure sheaves of reduced curves. This will be the only part of the argument needed in what follows.

Next we prove a lemma analogous to lemma \ref{conproj}:

\begin{lma}
\label{concurve}
Let $X$ be a noetherian variety and $\mcE$ a bounded complex of sheaves. Assume that
$$
\sum_{i=-\infty}^\infty \mathrm{dim\ }\Hom_{\mcD^b(\Coh(X))}(\mcO_C,\mcE[i]) < \infty
$$
for all $\mcO_C$ the structure sheaves of reduced curves in the support of $\mcE$. Then,
$$
\sum_{i=-\infty}^\infty \mathrm{dim\ }\Hom_{\Coh(X)}(\mcO_C,H^i(\mcE)) < \infty\ .
$$
\end{lma}

\begin{proof}

As in the proof of lemma \ref{conproj}, we have a long exact sequence:
\begin{equation*}
\begin{split}
0 \longrightarrow \Hom(\mcO_C,H^0(\mcE))& \longrightarrow \Hom(\mcO_C,\mcE) \longrightarrow \Hom(\mcO_C,\tau_{\ge 1}\mcE) \longrightarrow  \\
& \Hom(\mcO_C,H^0(\mcE)[1]) \longrightarrow\Hom(\mcO_C,\mcE[1]) \longrightarrow \Hom(\mcO_C, (\tau_{\ge 1}\mcE)[1]) \longrightarrow \dots \ .
\end{split}
\end{equation*}
From the beginning of this sequence, we have that $\Hom(\mcO_C,H^0(\mcE))$ is finite dimensional. By the argument in lemma \ref{supp}, $H^0(\mcE)$ has proper support and, hence, $\Hom(\mcO_C,H^0(\mcE)[i])$ is finite dimensional for all $i$. Thus, we again see that $\tau_{\ge 1}\mcE$ satisfies our hypothesis, and we proceed as in lemma \ref{conproj}.
\end{proof}

Finally, we dispatch with proposition \ref{sheaf}:
\setcounter{iprop}{\value{sheafc}}
\begin{iprop}
Let $X$ be a smooth noetherian variety. Then $\mcD^c_{cs}(\Coh(X))$ is equivalent to the full subcategory of $\mcD^c(\Coh(X))$ consisting of objects whose cohomology sheaves have proper support.
\end{iprop}
\begin{proof}
%Let $\mcE$ be a perfect complex of sheaves whose cohomology sheaves have proper support, and let $\mcF$ by any perfect complex. As above, there exists a spectral sequence converging to $\Hom(\mcF,\mcE[n])$ with
%$$
%E_2^{pq} =  \bigoplus_{i-j=p} \Ext_{\Coh(X)}^q(H^i(\mcF),H^j(\mcE))\ .
%$$
%As $X$ is smooth, we can represent $\mcE$ and $\mcF$ as bounded complexes of locally free sheaves. By lemma \ref{supp}, the spaces in the above $E_2$ term are all finite dimensional. Since $X$ is smooth and $\mcE$ and $\mcF$ are bounded, only a finite number are nonzero. Thus $\mcE$ has compact support.

Since $X$ is smooth, all perfect complexes are bounded. Thus, one direction follows immediately from proposition \ref{supsheaf}.

For the other direction, assume that $\mcE$ has compact support. Since $X$ is smooth, the structure sheaf of a reduced curve $C$ is a perfect object, and the hypotheses of lemma \ref{concurve} are satisfied. Thus, for all $C$ we have
$$
\mathrm{dim\ }\Hom_{\Coh(X)}(\mcO_C,H^i(\mcE)) < \infty\ .
$$
It then follows from the argument of lemma \ref{supp} that each sheaf $H^i(\mcE)$ has proper support, and we are done.
\end{proof}

\bibliographystyle{utphys}
\bibliography{thebib}

\end{document}